\documentclass[]{article}
\usepackage[utf8]{inputenc}

\usepackage[fleqn]{amsmath}
\usepackage{verbatim,amsthm}
\usepackage{graphicx}
\usepackage[colorinlistoftodos]{todonotes}
\usepackage[colorlinks=true, allcolors=blue]{hyperref}
\usepackage{amssymb,xcolor,amscd}
\usepackage{float}
\usepackage{enumitem}

\newcommand{\Rd}{\mathbb{R}}

\newcommand{\Un}{\mathbf{1}}

\DeclareMathOperator{\Prob}{P}
\DeclareMathOperator{\Mean}{E}

\let\0=\varnothing
\let\1=\Un

\theoremstyle{plain}
 \newtheorem{thm}{\bfseries Theorem.}[section]
 
 \newtheorem{prop}{\bfseries Proposition.}[section]
 
\theoremstyle{remark}

\numberwithin{equation}{section}

\usepackage{float}
\usepackage{listings}
\usepackage{enumitem}
\usepackage{authblk}

\title{On a Couple of Unresolved Group Testing Conjectures}

\author[1]{Ugnė Čižikovienė}
\author[1]{Viktor Skorniakov\thanks{corresponding author; e-mail: viktor.skorniakov@mif.vu.lt}}
\affil[1]{Institute of Applied Mathematics, Faculty of Mathematics and Informatics, Vilnius University, Naugarduko 24, Vilnius LT-03225, Lithuania}

\date{}

\begin{document}

\maketitle

\begin{abstract}
  In the recent review published in 2019, Malinovsky and Albert conjectured analytical formulae of the optimal sample sizes for the modified Dorfman and Sterret group testing schemes and verified the validity of the formulae numerically within a broad range. However, the authors pointed out that they were unable to give a rigorous proof of the stated conjectures. In this note, we affirm and refine their hypotheses analitically thereby filling in the present gap.
\end{abstract}

\section{Introduction}\label{s:intro}

Assume that $N$ objects are drawn at random from the population which is known to contain items of two types: the good ones and the bad ones with the prevalence of the latter being equal to $p\in(0,1)$. An aim is to identify the bad items. For this, a sample measurement from each object is obtained. Based on the test of that sample, an object can be unambiguously classified as being good or bad. 

There are many ways to arrange separation of bad and good items. One of widely known strategies is called Group Testing (GT). The essential idea of GT is to replace testing of single items by testing of groups having more than one item. A basic and most known approach was suggested by Dorfman \cite{dorfman_detection_1943} who was looking for the way to save testing resources devoted for the identification of syphilis in the military environment while screening the soldiers during World War II. He proposed the following scheme. Given $N$ blood samples, pool them and test the pooled sample for the presence of bacteria. In case pool tests negatively, do nothing; in case it tests positively, retest each single sample to identify those being ill. The possibility of saving stems from the intuitive reasoning that, having small prevalence $p$, the pooled sample will test negatively quite often resulting thus in one test per batch spanning $N$ items instead of $N$ tests applied for each individual directly. However, a natural question occurs: given fixed $p\in(0,1)$, what value of $N=N(p)$ one should choose in order to have maximal possible savings under the assumed model? 

In general, the question is usually addressed in the following way. Having fixed GT scheme (with the Dorfman's scheme above being one among many possible choices), let $T=T(N,p)$ denote the total (random) number of tests applied per batch spanning $N$ objects. Consider
\begin{equation}\label{e:t_per_item}
    t=t(N,p)=\frac{\Mean T}{N}=\text{"average number of test per one object"},
\end{equation}
and choose value of $N=N_{opt}$ minimizing $N\mapsto t(N,p)$. Such value is called the \emph{optimal configuration} and is one of the most important features characterising each particular GT scheme.

Turning to the case of the Dorfman scheme, $N_{opt}$ is known and was derived by Samuels \cite{samuels_exact_1978} who has shown that, for $p\in(0,1-(1/3)^{1/3})\approx (0,0.31)$, $N_{opt}\in\{\lfloor\sqrt{p^{-1}}\rfloor+1,\lfloor\sqrt{p^{-1}}\rfloor+2\}$  (here and in the sequel, $\lfloor x \rfloor$ stands for an integer part of $x\in\Rd$) whereas, for $p\geq 1- (1/3)^{1/3}$, $N_{opt}=1$. However, it appears that there is a couple of modifications of the Dorfman scheme for which analytical expression of $N_{opt}$, though specified explicitly as conjectures, left unverified rigorously (see \cite{malinovsky_revisiting_2019}). It is an aim of our note to fill in this gap. 

Finishing our short introduction, we note that, in addition to the \emph{probabilistic} GT described above, there is another subset of GT schemes altogether known under the name of \emph{combinatorial} GT schemes. Here, instead of assuming that the data was generated by some random mechanism, one assumes that the tested set of objects has $d$ bad items and tries to identify them by making use of combinatorial approach. Du and Hwang \cite{du_combinatorial_2000}, Ch.1, give a good introductory account on the history and development of both approaches. Our aim here is to emphasize that GT in general (and probabilistic GT in particular) is an important tool having many applications. The shortlist includes (yet is not limited to) the screening for other infectious diseases like HIV, hepatitis and COVID-19 (\cite{Xi-informativeness-1995}, \cite{wein_pooled_1996}, \cite{bilder_informative_2010}, \cite{may_pooled_2010}, \cite{stramer_nucleic_2011}, \cite{tebbs_two-stage_2013}, \cite{Shani_Narkiss_2020}, \cite{Sinnott-Armstrong_2020}), quality control (\cite{sobel+groll:1959}, \cite{Johnson1988StatisticalEO}, \cite{Johnson1990StatisticalEO}, \cite{bar-lev1990}, \cite{johnson_inspection_1991}), communication and security networking (\cite{wolf_born_1985}, \cite{laarhoven_efficient_2013}, \cite{skoric_generalised_2015}), genetics (\cite{du_pooling_2006}, \cite{Cutler41}, \cite{najafi2016fundamentalDNA}, \cite{Chang-chang-Cao16DNA}), and experimental physics (\cite{brady_searching_2000_exp_phys}, \cite{meinshausen_efficient_2009_exp_phys}, \cite{pletsch_optimized_2014_exp_phys}). 


The paper is organized as follows. In Section \ref{s:results}, we specify the afore mentioned modifications of the Dorfman scheme, shortly report on the work done in \cite{malinovsky_revisiting_2019} and other related references, and state our results. In Section \ref{s:conclusions}, some concluding remarks are provided. Finally, there are two appendixes placed in the end: Appendix \ref{a:proofs} contains proofs; Appendix \ref{a:figures} contains some accompanying figures.

\section{Task related background and results}\label{s:results}

\subsection{Notions and assumptions}
In what follows, we deal with two GT schemes: the modified Dorfman scheme and the Sterrett scheme (which is also a modification of the same Dorfman scheme). To facilitate exposition and distinguish between those two, we introduce abbreviations and write {D}, when referring to the first one, and {S}, when referring to the second one. Finally, ${\mathrm{D_0}}$ is reserved for the case of the original Dorfman's scheme. These abbreviations are also used as braced superscripts of the quantities considered when there is a need to emphasize the particular scheme. E.g., $t^{(D)}$, refers to $t$ given in \eqref{e:t_per_item} and corresponding to the case of the modified Dorfman's scheme.

All assumptions discussed in the Introduction remain valid, prevalence $p\in(0,1)$ is assumed to be fixed unless stated otherwise, and $q$ stands for $1-p$.

\subsection{The schemes}

\emph{Scheme {D}}. The original Dorfman scheme described in the introductory Section \ref{s:intro} obeys one obvious inconsistency. Given that pool needs retesting, it tests all objects even in the case when, having tested $N-1$ out of $N$, the first $N-1$ tested negatively. Sobel and Groll \cite{sobel+groll:1959} suggested not to test the last object in the latter case since his outcome is already predefined. Combining this modification with the description of the original Dorfman's scheme, it follows that, for $N\geq 2$,
\begin{equation}\label{e:tD_N}
    t^{(D)}(N,p)=\frac{\Mean T^{(D)}(N,p)}{N} = 1-q^N + \frac{1-pq^{N-1}}{N}
\end{equation}
since $\Prob(T^{(D)}=1)=q^{N}, \Prob(T^{(D)}=N)=pq^{N-1}$, and (consequently) $\Prob(T^{(D)}=N+1)=1-q^N-pq^{N-1}$.

\smallskip\noindent\emph{Scheme {S}}. Sterrett \cite{sterrett_detection_1957} suggested to modify {the scheme D} by retesting initial positive pool until appearance of the first positive case and then again apply pool testing to the remaining tail. In case the tail tests positive, one should proceed recursively as previously until, at some stage, tail tests negative or the whole set of individuals gets tested. 

Sobel and Groll \cite{sobel+groll:1959} were first to demonstrate that
\begin{equation}\label{e:tS_N}
    t^{(S)}(N,p)= 2 - q + \frac{2q-p^{-1}(1-q^{N+1})}{N}.
\end{equation}
For a short alternative proof see \cite{malinovsky_revisiting_2019}, Appendix A.

\subsection{Previous work, conjectures and our results}\label{ss:previous_and_res}

Before proceeding to the discussion of related work, we mention the fundamental Ungar's result given in \cite{Ungar-1960}. It states that, in the case of the binomial setting considered, there is no GT scheme that is better than individual one-by-one testing provided $p\geq \frac{3-\sqrt{5}}{2}$. This explains the role and frequent appearance of the cut-off point $\frac{3-\sqrt{5}}{2}$ in all what follows.

Malinovsky and Albert \cite{malinovsky_revisiting_2019} have demonstrated that, for a fixed $p\in(0,\frac{3-\sqrt{5}}{2})$, function $[1,\infty)\ni N\mapsto t^{(S)}(N,p)$ admits a unique absolute minimum which is attained at the unique zero of $[1,\infty)\ni N\mapsto\frac{\partial}{\partial N}t^{(S)}(N,p)$, say $N^{(S)}_*$. Therefore, $N_{opt}^{(S)}\in\{\lfloor N^{(S)}_*\rfloor,\lfloor N^{(S)}_*\rfloor+1\}$, and the choice between two possible values is made by evaluating whether $t^{(S)}({\lfloor N^{(S)}_*\rfloor},p)>t^{(S)}({\lfloor N^{(S)}_*\rfloor+1,p})$ or $t^{(S)}({\lfloor N^{(S)}_*\rfloor},p)\leq t^{(S)}({\lfloor N^{(S)}_*\rfloor+1},p)$ holds true.

For the case of the scheme {D}, Pfeifer and Enis \cite{pfeifer_dorfman-type_1978} have obtained a quite similar result stated below. 

\begin{thm}\label{t:pfeifer_enis}[Pfeifer and Enis \cite{pfeifer_dorfman-type_1978}, Lemma 2]
For a fixed $p\in\left(0,\frac{3-\sqrt{5}}{2}\right)$, function $[1,\infty)\ni N\mapsto t^{(D)}(N,p)$ admits a unique absolute minimum attained at the {smallest} zero of $[1,\infty)\ni N\mapsto\frac{\partial}{\partial N}t^{(D)}(N,p)$, denoted, in what follows, by $N^{(D)}_*$. In the set $A^{(D)}=\left\{(p,N)\in \left(0,\frac{3-\sqrt{5}}{2}\right)\times[1,\infty):t^{(D)}(N,p)<1\right\}$, $N^{(D)}_*$ is the only zero of $N\mapsto \frac{\partial}{\partial N}t^{(D)}(N,p)$.
\end{thm}
It is straightforward to check that $t^{(D_0)}(N,p)=1-q^N + \frac{1}{N}$. Therefore, by \eqref{e:tD_N}, $t^{(D_0)}(N,p)-t^{(D)}(N,p)=\frac{pq^{N-1}}{N}>0$. Hence, appealing to the result of Samuels \cite{samuels_exact_1978} mentioned in the Introduction, we thus conclude that looking for $N^{(D)}_{opt}$ corresponding to $p\in (0,1-(1/3)^{1/3})$, one can apply the same algorithm as for $N^{(S)}_{opt}$. Namely, it is enough to find the unique $N_*^{(D)}$ and then select $N^{(D)}_{opt}\in\{\lfloor N^{(D)}_*\rfloor,\lfloor N^{(D)}_*\rfloor+1\}$. In the region $\left[1-(1/3)^{1/3},\frac{3-\sqrt{5}}{2}\right)$ additional care is needed and we address that issue below.

Our previous discussion does not include results on analytical expressions for $N_{opt}^{(D)}, N_{opt}^{(S)}$. It appears that these are not known. Being unable to give analytical proof and based on numerical investigations, Malinovsky and Albert \cite{malinovsky_revisiting_2019} however conjectured that, for a fixed $p\in(0,\frac{3-\sqrt{5}}{2})$,
\begin{align}
    & N_{opt}^{(D)}\in\left\{\lfloor\sqrt{p^{-1}}\rfloor,\lfloor\sqrt{p^{-1}}\rfloor+1\right\}\quad\text{and}\\
    & N_{opt}^{(S)}\in\left\{\lfloor\sqrt{2p^{-1}}\rfloor,\lfloor\sqrt{2p^{-1}}\rfloor+1,\lfloor\sqrt{2p^{-1}}\rfloor+2\right\}.
\end{align}
Our results stated in two propositions below {affirm} these conjectures.
\begin{prop}\label{p:sterret}
Let 
\begin{equation}
    g_0(p):=\frac{1}{q}\left( \frac{1-2pq}{q\left(1-\ln q\sqrt{\frac{2}{p}}\right)} \right)^{\sqrt{\frac{p}{2}}}\quad \text{for}\quad p\in\left(0,\frac{3-\sqrt{5}}{2}\right).
\end{equation}
$g_0^{-1}{(\{1\})}$ consist of a single point $p_*$ an approximate value of which is 0.1711.
$N_{opt}^{(S)}\in\left\{{\lfloor\sqrt{2p^{-1}}\rfloor,\lfloor\sqrt{2p^{-1}}\rfloor+1}\right\}$ for {$p\in\left(p_*,\frac{3-\sqrt{5}}{2}\right)$}; 
    $N_{opt}^{(S)}\in\Big\{\lfloor\sqrt{2p^{-1}}$, $\lfloor\sqrt{2p^{-1}}\rfloor+1, \lfloor\sqrt{2p^{-1}}\rfloor+2\Big\}$ for {$p\in (0,p_*]$}. 
\end{prop}
\begin{prop}\label{p:dorfman}
For all $p\in\left(0,\frac{3-\sqrt{5}}{2}\right)$,
\begin{equation}\label{e:dorfman_opt}
    N_{opt}^{(D)}\in\left\{\lfloor\sqrt{p^{-1}}\rfloor,\lfloor\sqrt{p^{-1}}\rfloor+1\right\}.
\end{equation}
\end{prop}

\section{Concluding remarks}\label{s:conclusions}

The results available up to now allowed to obtain numerical solutions of optimal configurations $N_{opt}^{(S)}, N_{opt}^{(D)}$. Analytical expressions given in Propositions \ref{p:sterret}, \ref{p:dorfman} facilitate this task. This is one application. Another one, more important in our opinion, is an ability to compare analytically the schemes considered here against each other and to others not considered in terms of gain attained at optimal configuration. 

{To have a quick example, consider asymptotic regime when $p\to0+0$. It is then an exercise to show that, under this assumption, $t^{(D)}(N_{opt}^{(D)},p)=\frac{{2}+o(1)}{N_{opt}^{(D)}}$ and $t^{(S)}(N_{opt}^{(S)},p)=\frac{{2}+o(1)}{N_{opt}^{(S)}}$ as $p\to0+0$.} Therefore, taking into account Propositions \ref{p:sterret}--\ref{p:dorfman},
\begin{equation}\label{e:MA_limit}
    \lim_{p\to0+0} \frac{t^{(D)}(N_{opt}^{(D)},p)}{t^{(S)}(N_{opt}^{(S)},p)}=
    \lim_{p\to0+0}\frac{N_{opt}^{(S)}}{N_{opt}^{(D)}}=\sqrt{2}.
\end{equation}
Hence, for small $p$'s, testing optimally by making use of scheme D results in an average number of tests per object which is approximately 1.4 times larger than the average number of tests per object obtained applying scheme S at its optimal configuration. {Analyses of similar kind prevail in the literature. E.g., in \cite{malinovsky2021nested}, the limit \eqref{e:MA_limit} appears with D replaced by $\mathrm{D_0}$ and is the same. This means that asymptotically the original and the modified Dorfman schemes are identical. However, their behaviour is different for $p$'s far from the origin (see \cite{malinovsky_revisiting_2019}).}

Finishing, we note that, though developed quite long ago, schemes D and S are not outdated and still in use\footnote{e.g., the American Red Cross makes use of Dorfman scheme for the screening
of blood donations for HIV and hepatitis \cite{dodd_current_2002} {whereas, in Lithuania, Dorfman scheme is currently applied to test for COVID-19 employees of larger firms and pupils attending public schools \cite{PoolLT}}} even in large scale projects. A reason is quite simple. Despite the fact that since original work of Dorfman many schemes tied to particular needs of applications considered were developed, the operational characteristics of the test kit and organisational flow of the whole project may not afford to apply more elaborated schemes. In such cases, simple schemes appear to be a good alternative resulting in cost savings.

\appendix
\section{Proofs}\label{a:proofs}
\begin{figure}[ht!]
\centering
\includegraphics[width=15cm, height=12cm]{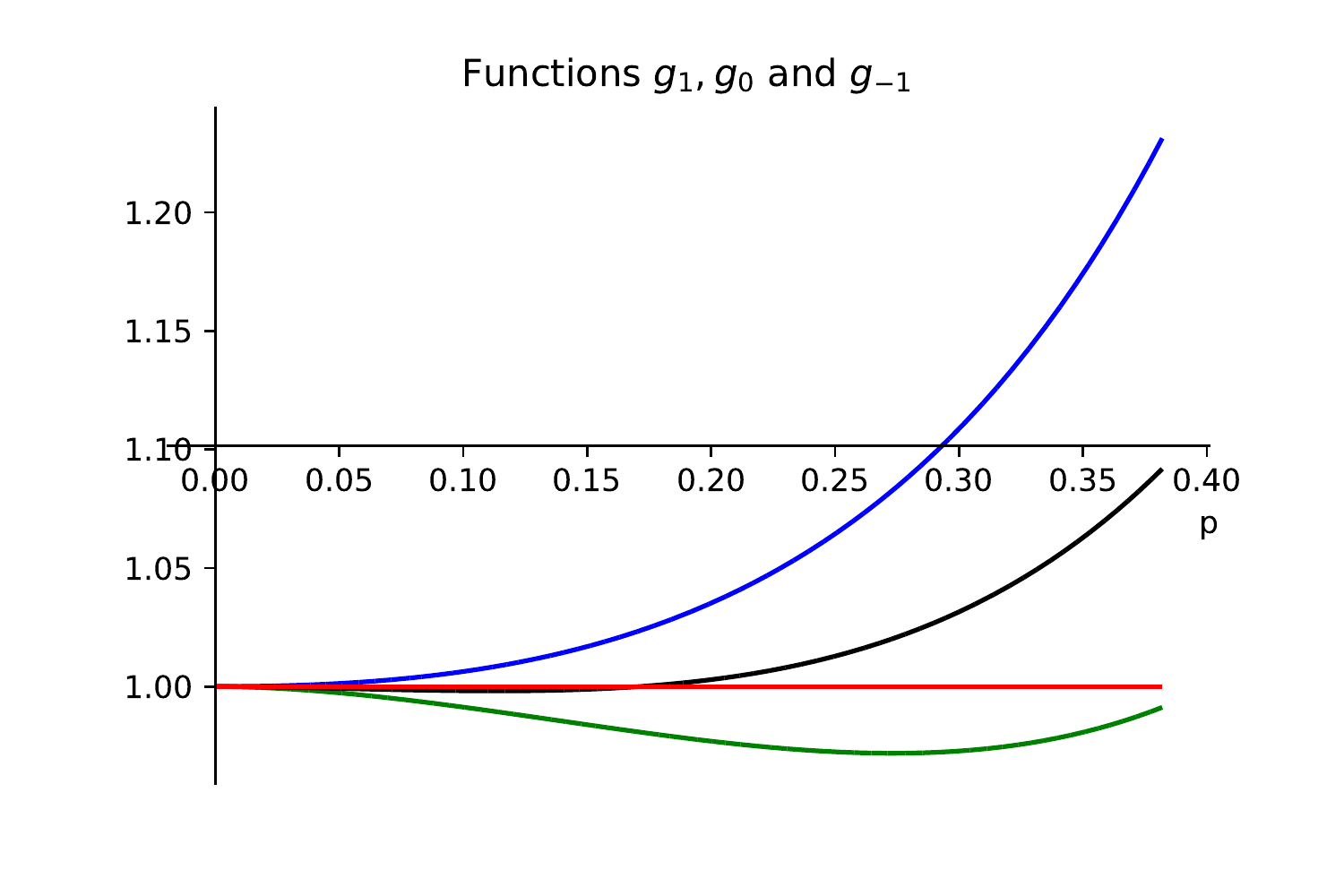}
\caption{Functions $g_m$ on $\left(0,\frac{3-\sqrt{5}}{2}\right)$: $g_1$ plotted in blue, $g_0$ plotted in black, and $g_{-1}$ plotted in green.} \label{fig:all_g}
\end{figure}

\begin{figure}[ht!]
\centering
\includegraphics[width=15cm, height=12cm]{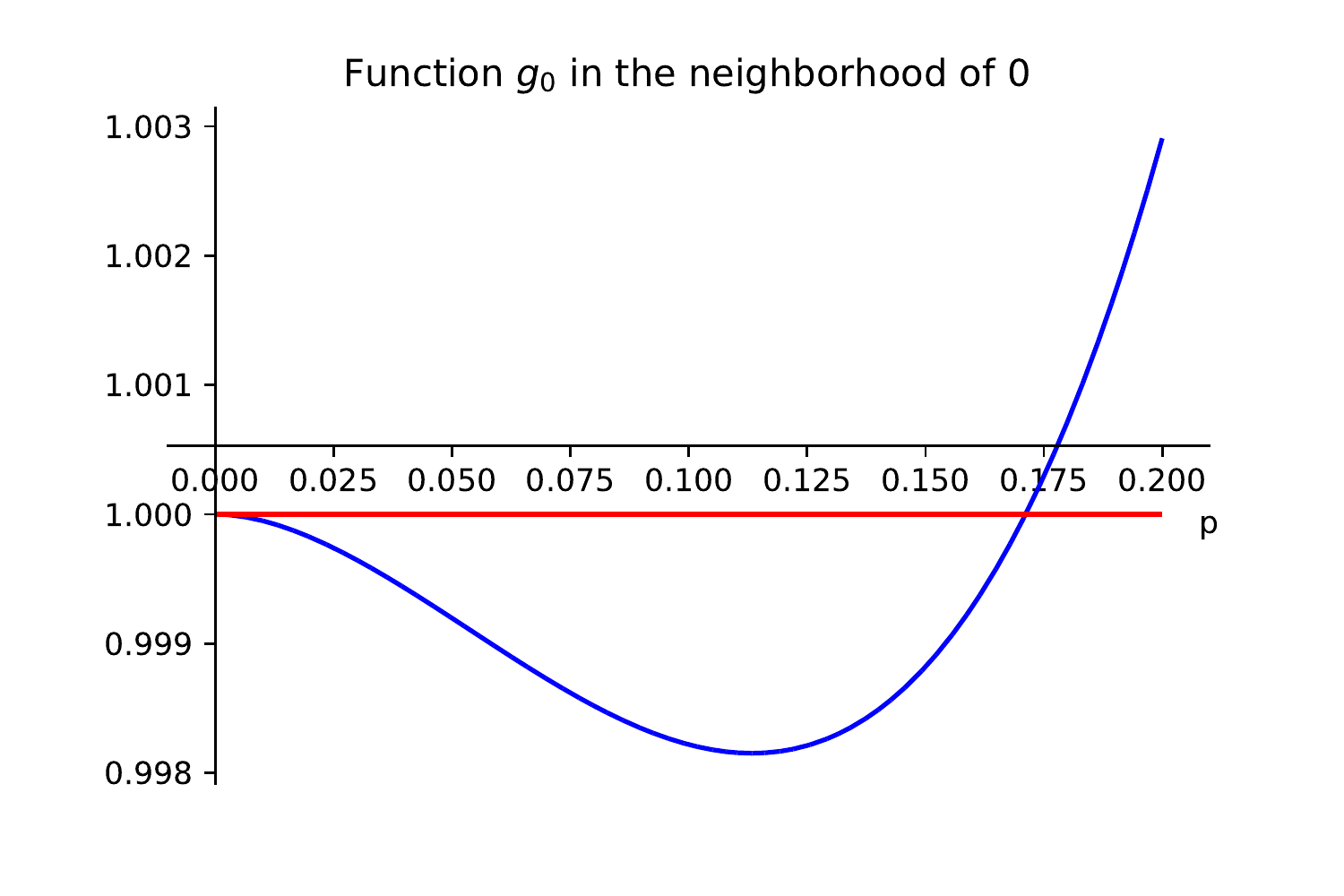}
\caption{The behaviour of $g_0$ near the origin: the point of intersection with 1 is the cut-off point $p_*$. For $p$'s on the left from it, $N_{opt}^{(S)}\in\{{\lfloor\sqrt{2p^{-1}}\rfloor,\lfloor\sqrt{2p^{-1}}\rfloor+1,\lfloor\sqrt{2p^{-1}}\rfloor+2}\}$, whereas for $p$'s on the right, $N_{opt}^{(S)}\in\{{\lfloor\sqrt{2p^{-1}}\rfloor,\lfloor\sqrt{2p^{-1}}\rfloor+1}\}$} \label{fig:g_0}
\end{figure}

\emph{Proof of Proposition \ref{p:sterret}}. Consider equation $\frac{\partial}{\partial N}t^{(S)}(N,p)=0$. Simple rearrangement shows that it is equivalent to equality
\begin{equation}\label{e:sterret_fixed_point}
    \frac{1}{\ln q} - \frac{1-2pq}{\ln q}\left(\frac{1}{q}\right)^{N+1}=N.
\end{equation}
Denote the lhs by $h(N)$. Then \eqref{e:sterret_fixed_point} means that $N_*^{(S)}$ is a fixed point of $h:[0,\infty)\to\Rd$. Since $h$ is translated and scaled increasing exponential function with $h(0)<0$, that fixed point is unique in agreement with results discussed in Subsection \ref{ss:previous_and_res}. Moreover, in a view of these results, it suffices to demonstrate that $N_*^{(S)}\in\left[\sqrt{2p^{-1}}-1,\sqrt{2p^{-1}}+1\right]$ in order to deduce that $N^{(S)}_{opt}\in\left\{\lfloor\sqrt{2p^{-1}}\rfloor+i:i\in\{-1,0,1,{2}\}\right\}$. Taking into account the exponential form of $h$, the latter will follow provided we show that
\begin{equation}\label{e:sterret_main_ineq}
    h\left(\sqrt{\frac{2}{p}}+1\right)> \sqrt{\frac{2}{p}}+1\quad\text{and}\quad
    h\left(\sqrt{\frac{2}{p}}-1\right)< \sqrt{\frac{2}{p}}-1.
\end{equation}
For each $m\in\{-1,0,1\}$, define a function
\begin{equation}
    g_{m}(p)=\frac{1}{q}\left( \frac{1-2pq}{q^{1+m}\left(1-\ln q\sqrt{\frac{2}{p}}\left(1+m\sqrt{\frac{p}{2}}\right)\right)} \right)^{\sqrt{\frac{p}{2}}}\quad \text{for}\quad p\in\left(0,\frac{3-\sqrt{5}}{2}\right).
\end{equation}
By simple rearrangement, it follows that \eqref{e:sterret_main_ineq} is equivalent to
\begin{equation}\label{e:sterret_g_ineq}
    g_{1}(p)> 1\quad \text{and} \quad g_{-1}(p)< 1.
\end{equation}
Figure \ref{fig:all_g} shows the graphs of $g_{1}, g_0$, and $ g_{-1}$. These suggest that relationships \eqref{e:sterret_g_ineq} do hold outside the zero neighborhood though, due to resolution issues, the true behaviour close to origin may be masked. To see that here $g_{1}(p)>1$ and $g_{-1}(p)<1$, consider expansion
\begin{multline}
    \ln g_{m}(p)=
    \sqrt{\frac{p}{2}}\left(
      \ln(1-2pq)-(1+m)\ln q - \ln\left(1-\ln q\left(m+\sqrt{\frac{2}{p}}\right) \right)
    \right)\\ - \ln q=
    \sqrt{\frac{p}{2}}\left(
      p(1+m-2q) + p^2\left(\frac{1+m}{2}-2q^2\right) + O(p^3) 
    \right) + \\
    m\sqrt{\frac{p}{2}}\ln q + \sqrt{\frac{1}{2p}}(\ln q)^2\left(1+m\sqrt{\frac{p}{2}}\right)^2 + 
    \frac{2}{3}\frac{(\ln q)^3}{p}\left(1+m\sqrt{\frac{p}{2}}\right)^3 + \\
    \frac{1}{\sqrt{2}}\frac{(\ln q)^4}{p^{3/2}}\left(1+m\sqrt{\frac{p}{2}}\right)^4 +
    O\left(p^3\right).
\end{multline}
Plugging in $m=\pm 1$ yields
\begin{multline*}
    \ln g_{-1}(p)=
    \sqrt{\frac{p}{2}}\left(
     -2pq(1+ pq + O(p^2)
    \right) - 
    \sqrt{\frac{p}{2}}\ln q\left(1 - {\frac{\ln q}{p}}\right) - (\ln q)^2 +\\ \frac{\sqrt{p}(\ln q)^2}{2\sqrt{2}} + 
    \frac{2}{3}\frac{(\ln q)^3}{p}\left(1-\sqrt{\frac{p}{2}}\right)^3 + 
    \frac{1}{\sqrt{2}}\frac{(\ln q)^4}{p^{3/2}}\left(1-\sqrt{\frac{p}{2}}\right)^4 +
    O\left(p^3\right)=\\
    -\sqrt{\frac{p}{2}}\left(
        2p + \ln q\left(2+O(p)\right)
    \right) - 
    (\ln q)^{2}\left(
        1 - \frac{2}{3}\frac{\ln q}{p}\left(1-\sqrt{\frac{p}{2}}\right)^3
    \right)+ O\left(p^{5/2}\right)\\=
    -\frac{5}{3}p^2 + O\left(p^{5/2}\right).
\end{multline*}
and
\begin{multline*}
    \ln g_1(p)=
    \sqrt{\frac{p}{2}}\left(
     3p^2-2(pq)^2 + O(p^3)
    \right) + \sqrt{\frac{p}{2}}\ln q\left(1 + {\frac{\ln q}{p}}\right) + (\ln q)^2 +\\ \frac{\sqrt{p}(\ln q)^2}{2\sqrt{2}} + 
    \frac{2}{3}\frac{(\ln q)^3}{p}\left(1+\sqrt{\frac{p}{2}}\right)^3 + 
    \frac{1}{\sqrt{2}}\frac{(\ln q)^4}{p^{3/2}}\left(1+\sqrt{\frac{p}{2}}\right)^4 +
    O\left(p^3\right)=\\
    (\ln q)^{2}\left(
        1 + \frac{2}{3}\frac{\ln q}{p}\left(1+\sqrt{\frac{p}{2}}\right)^3
    \right)=\frac{1}{3}p^2 + O\left(p^{5/2}\right).
\end{multline*}
Hence, $g_1(0+0)=g_{-1}(0+0)=1$ and $g_1(p)>1,g_{-1}(p)<1$ for all $p\in(0,\delta)$ provided $\delta>0$ is small enough. One can further show that $g_1^{\prime}$ is positive on $(0,\frac{3-\sqrt{5}}{2})$, whereas, in case of $g_{-1}^{\prime}$, the following hold true: it has a unique zero $x_0\in(0,\frac{3-\sqrt{5}}{2})$; it is negative on $(0,x_0)$ and positive on $(x_0,\frac{3-\sqrt{5}}{2})$; finally, $\lim\limits_{x\to \frac{3-\sqrt{5}}{2}-0}g_{-1}\left(x\right)\approx 0.9912$. Calculations being lengthy and tedious nonetheless require only standard calculus and we, therefore, omit the details. Putting all together, relationships \eqref{e:sterret_g_ineq} do hold. 
Taking into account all the said and then by the similar argument as above, it follows that
\begin{gather*}
    N_*^{(S)}\in\left[\sqrt{2p^{-1}}-1,\sqrt{2p^{-1}}\right] \Longleftrightarrow g_{0}(p)\geq 1\\
    \quad \text{and} \quad \\
    N_*^{(S)}\in\left[\sqrt{2p^{-1}},\sqrt{2p^{-1}}+1\right] \Longleftrightarrow g_{0}(p)\leq 1.
\end{gather*}
Also, since $\left(0,\frac{3-\sqrt{5}}{2}\right)\ni p\mapsto N_*^{(S)}(p)$ is continuous and strictly decreasing\footnote{this needs some reasoning yet we omit the details}, $g_0-1$ has a unique zero $p_*\in(0,\frac{3-\sqrt{5}}{2})$ (see figure \ref{fig:g_0}) and sets $g_0^{-1}((-\infty,1]),g_0^{-1}([1,\infty))$ are connected, i.e., intervals $(0,g_0^{-1}(\{1\})],\left[g_0^{-1}(\{1\}),\frac{3-\sqrt{5}}{2}\right)$ respectively. 

\begin{figure}[ht!]
\centering
\includegraphics[width=15cm, height=12cm]{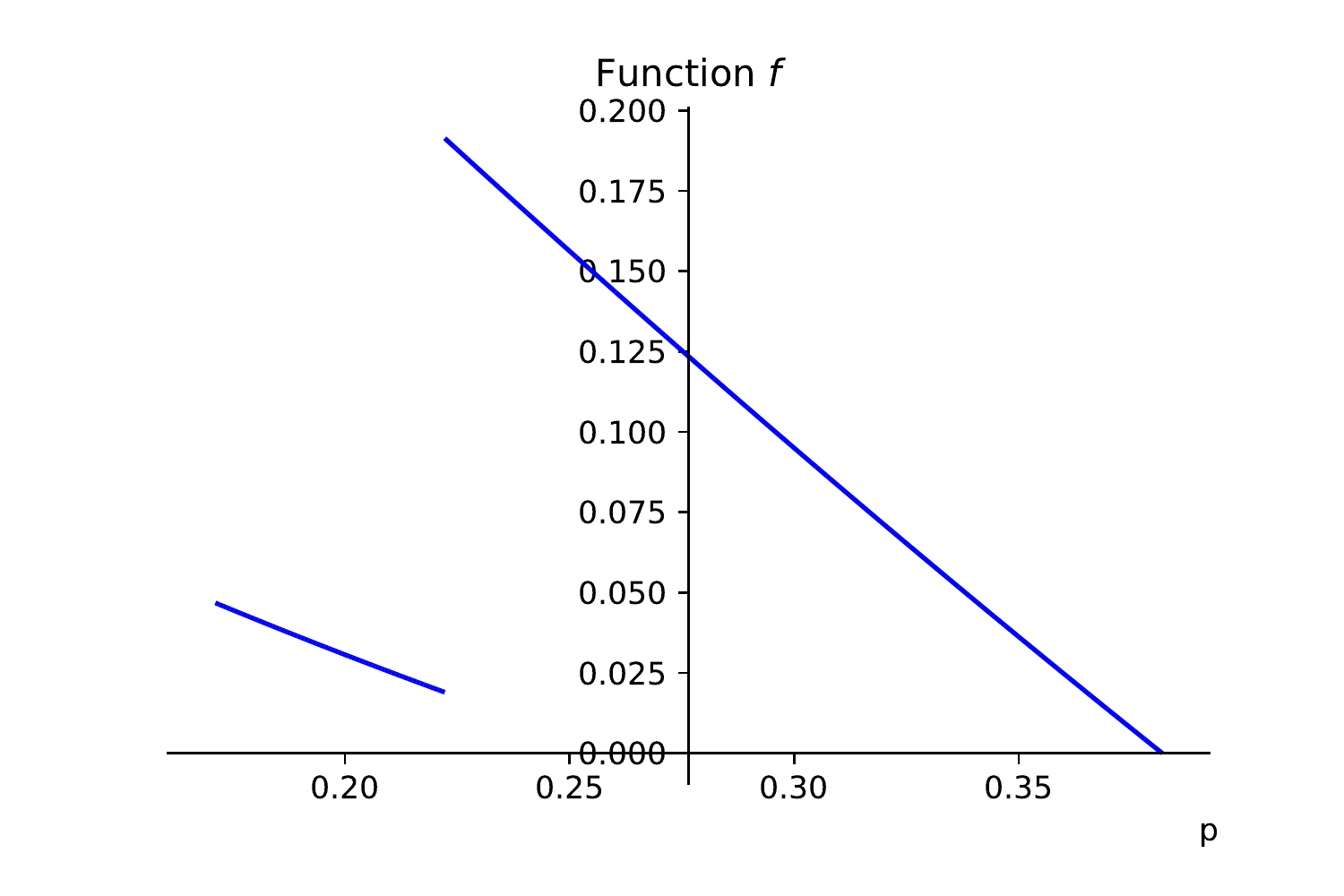}
\caption{The behaviour of function $f(p)=t^{(S)}\left(\left\lfloor\sqrt{\frac{2}{p}}\right\rfloor-1,p\right)-
    t^{(S)}\left(\left\lfloor\sqrt{\frac{2}{p}}\right\rfloor,p\right)$ on $\left(p_*,\frac{3-\sqrt{5}}{2}\right)$.} 
\label{fig:f}
\end{figure}

{To finish the proof, it remains to demonstrate that $\left\lfloor\sqrt{\frac{2}{p}}\right\rfloor-1$ is never optimal on $\left(p_*,\frac{3-\sqrt{5}}{2}\right)$. For this, consider function
\begin{multline*}
    f(p)=t^{(S)}\left(\left\lfloor\sqrt{\frac{2}{p}}\right\rfloor-1,p\right)-
    t^{(S)}\left(\left\lfloor\sqrt{\frac{2}{p}}\right\rfloor,p\right)=\\
    \1_{\left(p_*,\frac{2}{9}\right]}(p)\left(t^{(S)}\left(1,p\right)-
    t^{(S)}\left(2,p\right)\right)+
    \1_{\left(\frac{2}{9},\frac{3-\sqrt{5}}{2}\right)}(p)\left(t^{(S)}\left(2,p\right)-
    t^{(S)}\left(3,p\right)\right)
\end{multline*}
on $\left(p_*,\frac{3-\sqrt{5}}{2}\right)$. It is left continuous and has a single point of discontinuity equal to $\frac{2}{9}$. Since $f^\prime$ is negative on $\left(p_*,\frac{2}{9}\right)\cup \left(\frac{2}{9},\frac{3-\sqrt{5}}{2}\right)$, invoking left continuity, we infer that its minimal values
are $f(2/9)$ on $(p_*,2/9]$ and $f\left(\frac{3-\sqrt{5}}{2}-0\right)$ on $\left(\frac{2}{9},\frac{3-\sqrt{5}}{2}\right)$ respectively. By direct substitution and because of left continuity, $f\left(\frac{3-\sqrt{5}}{2}-0\right)=0$, whereas numerical estimation yields $f(2/9)\approx 0.018976$. The verification of the enumerated properties of $f$ requires only lengthy standard calculus and we omit it. The graph of $f$ illustrating its behaviour is given in figure \ref{fig:f}.}
\qed

\medskip\emph{Proof of Proposition \ref{p:dorfman}}. For the sake of clarity, we split the proof into three steps. As in the proof of Prop. \ref{p:sterret}, some tedious details are omitted and only the sketch is given.

\smallskip\emph{Step 1: slicing $A^{(D)}$.} In this step, we show that
\begin{equation}\label{e:bound_for_tN_less_1}
    \left\{\left( p, \frac{1}{\sqrt{p}} + 1 - \frac{5}{2}p \right):p\in\left(0,\frac{3-\sqrt{5}}{2}\right)\right\}\subseteq
    A^{(D)},
\end{equation}
where $A^{(D)}$ is the same as in the statement of Theorem \ref{t:pfeifer_enis}. To this end, note that, by elementary rearrangement,
\begin{multline*}
    t^{(D)}(N,p)<1 \Longleftrightarrow  1< q^{N-1}(Nq+p)=q^{N-1}((N-1)q+1)=\\
    \left[\text{denoting } N-1=x\right]=q^x(1+qx)
    \Longleftrightarrow 0<x\ln q + \ln(1+qx).
\end{multline*}
Plugging in $N=\sqrt{p^{-1}}+1-(5/2)p$, we then obtain the condition we need to check:
\begin{equation*}
    \left(\frac{1}{\sqrt{p}}-\frac{5p}{2}\right)\ln q+\ln\left(1+q\left(\frac{1}{\sqrt{p}}-\frac{5p}{2}\right)\right)>0.
\end{equation*}
Putting $y=\sqrt{p}$ this translates to checking that 
\begin{equation*}
    f(y)=\left(\frac{1}{y}-\frac{5y^2}{2}\right)\ln(1-y^2)+
    \ln\left(1+\left(1-y^2\right)\left(\frac{1}{y}-\frac{5y^2}{2}\right)\right)
\end{equation*}
is positive on its domain $\left(0,\sqrt{\frac{3-\sqrt{5}}{2}}\right)$. One can show that $f$ is strictly convex on $\left(0,\sqrt{\frac{3-\sqrt{5}}{2}}\right)$. Consequently, $f^{\prime}$ is non-decreasing and upper bounded by $\lim\limits_{x\to \sqrt{\frac{3-\sqrt{5}}{2}}-0} f^{\prime}(x)\approx -1.66$. This, in turn, yields that $f$ is decreasing and lower bounded by $\lim\limits_{x\to \sqrt{\frac{3-\sqrt{5}}{2}}-0} f(x)\approx 0.024$.

\smallskip\emph{Step 2: bracing the optimal points.} Rewrite equation $\frac{\partial }{\partial N}t^{(D)}(N,p) = 0$ as follows:
\begin{equation*}
     -\frac{q}{p}\left( N^2+\frac{1}{\ln q}\left(\frac{1}{q}\right)^N\right)+\frac{1}{\ln q} = N,
\end{equation*}
and denote the function on the rhs by $f(N,p)$. Then each fixed point of $N\mapsto f(N,p)$ is the zero of $\frac{\partial }{\partial N}t^{(D)}(N,p)$. In particular, the statement applies to $N_*^{(D)}$. In this step, by making use of this observation, we show that $N_*^{(D)}$ (the detailed explanation of the latter fact is given in \emph{Step 3}) does not deviate a lot from $\sqrt{p^{-1}}$. To achieve the goal, we consider points 
\begin{multline}\label{e:N_theta}
    N=N(\theta)=\left(\frac{1}{\sqrt{p}}+1-\frac{5}{2}p\right)\theta + 
    \left(\frac{1}{\sqrt{p}}-p\right)(1-\theta)=\\
    \frac{1}{\sqrt{p}}-p+\theta\left(1-\frac{3}{2}p\right),\theta\in[0,1],
\end{multline}
and demonstrate that 
\begin{multline}\label{e:step2_res}
    \forall p\in\left(0,\frac{3-\sqrt{5}}{2}\right)\exists !\theta\in[0,1]: f(N(\theta),p)=N(\theta)\Longleftrightarrow\\
    \forall p\in\left(0,\frac{3-\sqrt{5}}{2}\right)\exists !\theta\in[0,1]: 
    p^{\frac{3}{2}}(f(N(\theta),p)-N(\theta))=0.
\end{multline}
For convenience, put $h(\theta,p)=p^{\frac{3}{2}}(f(N(\theta),p)-N(\theta))$. Calculating derivative yields
\begin{multline}
    \frac{\partial}{\partial \theta}h(\theta,p)=p^{\frac{3}{2}}\left(
    \frac{\partial}{\partial N}f(N,p)-1\right)\frac{\partial}{\partial \theta}N(\theta)=\\
    p^{\frac{3}{2}}\left(-\frac{q}{p}\left(2N(\theta)-\left(\frac{1}{q}\right)^{N(\theta)}\right) -1\right)\left(1-\frac{3}{2}p\right)
\end{multline}
and then, since $\frac{\partial}{\partial \theta}N(\theta)$ does not depend on $\theta$,
\begin{multline}\label{e:h_diff_theta2}
    \frac{\partial^2}{\partial \theta^2}h(\theta,p)=p^{\frac{3}{2}}\left(
    \frac{\partial^2}{\partial N^2}f(N,p)\right)\left(\frac{\partial}{\partial \theta}N(\theta)\right)^2=\\
    -{q}\sqrt{p}\left(2+\left(\frac{1}{q}\right)^{N(\theta)}\ln q \right)\left(1-\frac{3}{2}p\right)^2.
\end{multline}
For a fixed $p$, consider the term $\left(2+\left(\frac{1}{q}\right)^{N(\theta)}\ln q \right)$. Since $N(\theta)\uparrow$, it is upper bounded by $\left(2+\left(\frac{1}{q}\right)^{N(1)}\ln q \right)$. By making use of the second derivative test, one can show that $p\mapsto \left(2+\left(\frac{1}{q}\right)^{N(1)}\ln q \right)$ is strictly concave and has negative derivative on $\left(0,\frac{3-\sqrt{5}}{2}\right)$. Therefore, it is lower bounded by its left limit at $\frac{3-\sqrt{5}}{2}$. The value of the latter is approximately equal to $0.93$. It follows then from \eqref{e:h_diff_theta2} that, for any fixed $p\in \left(0,\frac{3-\sqrt{5}}{2}\right)$, $\theta\mapsto \frac{\partial^2}{\partial \theta^2}h(\theta,p)$ is negative on $[0,1]$. Hence, for any fixed $p\in \left(0,\frac{3-\sqrt{5}}{2}\right)$, $\theta\mapsto h(\theta,p)$ is concave and has therefore a decreasing derivative. Inspection of $p\mapsto \frac{\partial}{\partial \theta}h(0,p)$ reveals that $\theta\mapsto \frac{\partial}{\partial \theta}h(\theta,p)$ is negative for any fixed $p$ meaning that the range of $\theta\mapsto h(\theta,p)$ is equal to $[h(1,p),h(0,p)]$. Finally, omitting the details of a tedious exercise of verification that the range of $p\mapsto h(1,p)h(0,p)$ lies in $(-\infty,0)$, we finish proof of this step and conclude that \eqref{e:step2_res} indeed holds.

\begin{figure}[ht!]
\centering
\includegraphics[width=15cm, height=12cm]{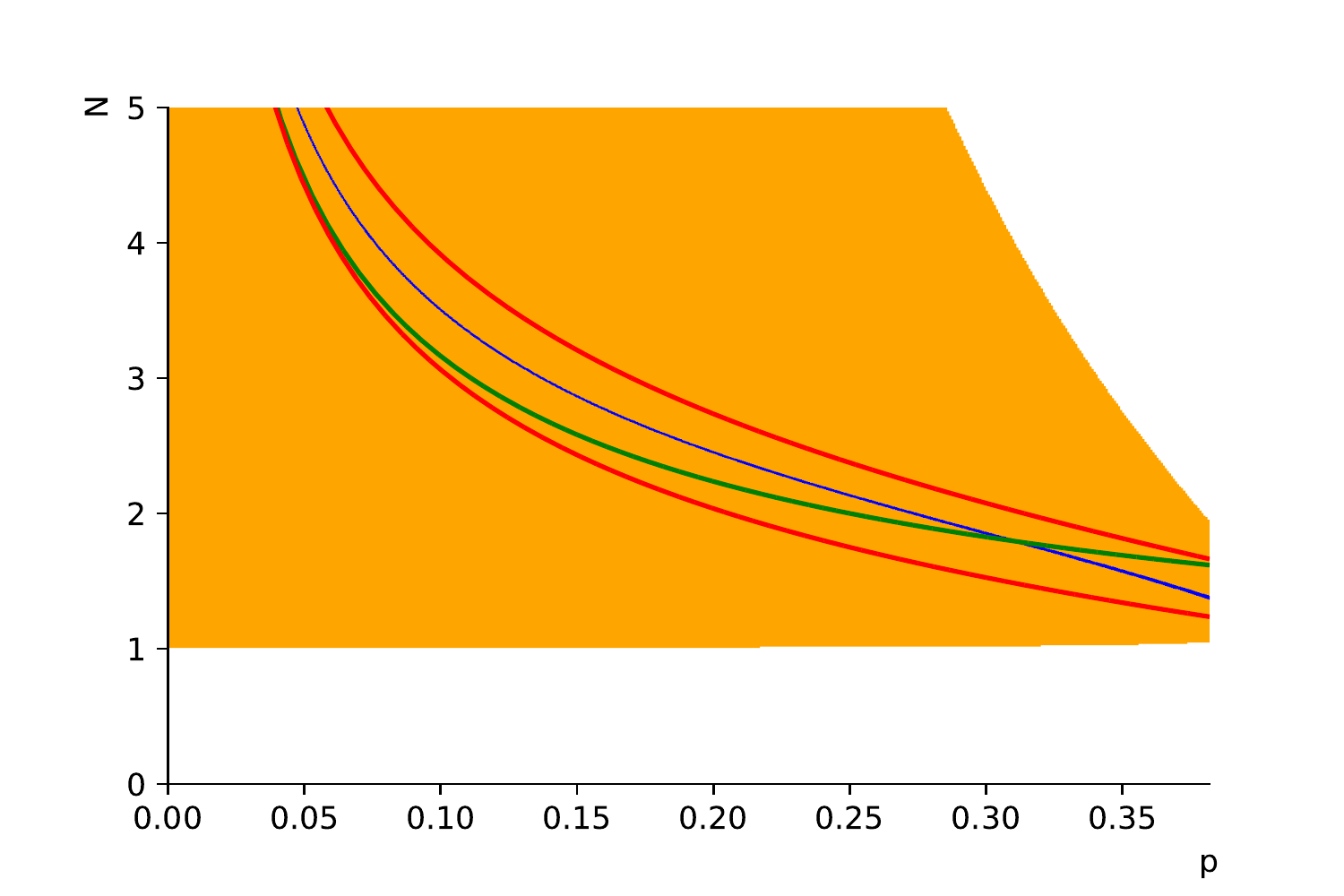}
\caption{orange area corresponds to $\{(N,p):t^{(D)}(N,p)<1\}$; upper and lower red lines show bracing functions $p\mapsto\sqrt{p^{-1}}+(1-(5/2)p)$ and $p\mapsto\sqrt{p^{-1}}-p$ respectively; green line shows function $p\mapsto\sqrt{p^{-1}}$; blue line shows $p\mapsto N^{(D)}_*$} \label{fig:Dorfman}
\end{figure}

\smallskip\emph{Step 3: the end of the proof.} By \emph{Step 2}, for each $p\in\left(0,\frac{3-\sqrt{5}}{2}\right)$, there exists $N(p)\in \left(\frac{1}{\sqrt{p}}-p,\frac{1}{\sqrt{p}}+1-\frac{5}{2}p\right)$ which solves $\frac{\partial }{\partial N}t^{(D)}(N,p) = 0$. By \emph{Step 1}, $N(p)\in A^{(D)}$. Therefore, by Theorem \ref{t:pfeifer_enis}, $N(p)$ is the unique global minimizer of $[1,\infty)\ni N\mapsto t^{(D)}(N,p)$, i.e., $N(p)=N^{(D)}_*$. This implies that $N_{opt}^{(D)}\in\{\lfloor\sqrt{p^{-1}}\rfloor-1,\lfloor\sqrt{p^{-1}}\rfloor,\lfloor\sqrt{p^{-1}}\rfloor+1\}$, and it remains to exclude the point $\lfloor\sqrt{p^{-1}}\rfloor-1$. The route is as follows. First, by making exactly the same technique as in \emph{Step 2}, show that $N(p)\in \left(\frac{1}{\sqrt{p}},\frac{1}{\sqrt{p}}+1-\frac{5}{2}p\right)$ for $p\in(0,0.3)$. Second, note that integer parts of $\frac{1}{\sqrt{p}}-p$ and $\frac{1}{\sqrt{p}}$ coincide for $p\in\left[0.3,\frac{3-\sqrt{5}}{2}\right)$. Figure \ref{fig:Dorfman} provides a good visual summary of the whole proof and explains the need of this workaround in particular. \qed

\section{Figures}\label{a:figures}
The figures below were produced by making use of an open source computer algebra system SymPy \cite{sympy}.

\clearpage
\bibliographystyle{alpha}
\newcommand{\etalchar}[1]{$^{#1}$}

\end{document}